\documentclass{article}

\usepackage[utf8]{inputenc}

\usepackage{amsfonts, amsmath, amssymb, amsthm}
\usepackage{graphicx}
\usepackage{comment}
\usepackage{url}
\usepackage{subfig}

\usepackage{multirow}
\usepackage{multicol}
\usepackage{booktabs}
\usepackage[table,xcdraw]{xcolor}

\usepackage{cite}

\begin{document}
\title{Basins of Attraction, Commitment Sets and Phenotypes of Boolean Networks}

\author{Hannes Klarner, Frederike Heinitz, Sarah Nee, Heike Siebert\\
Department of Mathematics, Freie Universität Berlin, Germany\\
E-mail: Hannes.Klarner@FU-Berlin.de}

\maketitle

\begin{abstract}
The attractors of Boolean networks and their basins have been shown to be highly relevant for model validation and predictive modelling, e.g., in systems biology. 
Yet there are currently very few tools available that are able to compute and visualise not only attractors but also their basins.
In the realm of asynchronous, non-deterministic modeling not only is the repertoire of software even more limited, but also the formal notions for basins of attraction are often lacking.
In this setting, the difficulty both for theory and computation arises from the fact that states may be elements of several distinct basins.
In this paper we address this topic by partitioning the state space into sets that are committed to the same attractors.
These commitment sets can easily be generalised to sets that are equivalent w.r.t. the long-term behaviours of pre-selected nodes which leads us to the notions of markers and phenotypes which we illustrate in a case study on bladder tumorigenesis.
For every concept we propose equivalent CTL model checking queries and an extension of the state of the art model checking software NuSMV is made available that is capable of computing the respective sets.
All notions are fully integrated as three new modules in our Python package PyBoolNet, including functions for visualising the basins, commitment sets and phenotypes as quotient graphs and pie charts.
\end{abstract}

\section{Background}\label{sec:introduction}

Boolean networks offer an intuitive approach to simulating the dynamics of interaction networks.
In cell biology these are usually gene regulatory or signal transduction networks.
Each biological component is modelled by a Boolean variable that can only switch between two states, \emph{true} and \emph{false}.
Depending on the component, these states represent whether a protein is present at high or low concentration levels,
a gene is being transcribed at or above its base rate, a signalling molecule is phosphorylated or not, and so on.
The conditions under which variables switch on and off are specified by Boolean functions that involve other variables of the network.
Together with a transition relation, which determines how many and which components are allowed to change during a state transition,
Boolean networks can be used to predict transient and long term activity profiles of the involved components.
Although the restriction to discrete states is arguably an over-simplification for many processes,
Boolean networks have frequently been shown to be fruitful to predict drug targets, phenotypes, novel interactions, etc.
This paper is particularly relevant to models whose state transitions are non-deterministic.
A prominent example are the so-called Thomas networks (after Ren\'e Thomas, see e.g. \cite{Thomas1991RegulatoryNetworks}) where changes in state variables are modelled as processes that are not synchronised.
Asynchronous transition systems are especially applicable when details about the precise mechanics of the processes, e.g. the activation and inhibition delays, are unknown or uncertain.
The results of this paper are particularly relevant to non-deterministic transition relations but also work in the deterministic case.
For a detailed introduction to Boolean networks and qualitative modelling more generally see \cite{Naldi2015Colomoto}.

Attractor detection is a frequent starting point for testing the consistency of models with data and also for obtaining model predictions, see e.g. \cite[Sec. 3.2]{WAJaoude2015ThPlasticity}.
An extension of this approach considers not only the location of attractors but also their basins of attraction.
These regions and transitions between them determine the likelihood of reaching particular attractors and not others.
From a practical, wet-lab point of view, the motivation for studying basins and attractor commitment sets is to find out how the underlying decision processes are implemented within the interaction networks of living cells.
Computational models and the formalisation of notions like basin, commitment set and phenotype add to this research by making hypothesis testable in silico.
When these models are validated they may then be used to measure the robustness of interaction networks with respect to knock-outs or knock-ins and other perturbations.
Due to the high computational cost of analysing the complete state transition graph for larger networks, such studies often use aggregated results from simulations  \cite{remy2015tumorigenesis}. 
While analytic results are mostly focused on Boolean networks with deterministic updates \cite{Willadsen2007StateSpaceStructure, Wuensche2010Basins},
there are still some algorithms for computing the basins of Boolean networks with asynchronous updates \cite{Bahi2006Basins,dinwoodie2016basins,Arellano2011Antelope}.
These all have to deal with the inherent difficulties in analysing asynchronous state transition graphs, which in turn results in restricted scalability.

In general, the computation of basins of attraction is challenging because the size of the state space grows exponentially with the number of variables of the model.
A competitive approach for dealing with large transition systems is to use so-called symbolic model checking algorithms \cite{Baier2008Principles}.
Model checking is a method from computer science to decide whether formal hardware designs satisfy desired specifications.
The designs are formalised as so-called Kripke structures, which are labelled, directed graphs,
and the specifications to be verified are translated into statements of a chosen temporal logic.
Computation Tree Logic (CTL) is an example of such a logic.
It allows statements about the alternative futures originating in states \cite{Baier2008Principles}.
These are particularly useful when dealing with properties related to attractors.
NuSMV is an open source symbolic model checking tool allowing for the efficient analysis of transition systems \cite{Cimatti2002NuSMV2}.

It has been observed that the state transition graphs (STGs) of Boolean networks are Kripke structures and that model checking is therefore applicable to Boolean networks \cite{Bernot2004TemporalLogic}.
There are custom model checking tools that are specifically developed for biological networks, e.g. \cite{streck2016toolkit,Arellano2011Antelope},
but their additional capabilities may come at the cost of efficiency in more basic queries.
Still, we believe that for the problem of computing basins of attractors and related objects state of the art software designed to answer reachability queries exactly such as the model checking software NuSMV provides the ideal basis for implementation.

In this paper we, on the one hand, formally develop several notions related to basins of attraction and, on the other hand, provide efficient tools to compute these objects in practice using a model checking approach.
Namely, we created NuSMV-a, an extension of NuSMV 2.6, which implements model checking with accepting states.
We complement the theoretical part in later sections by giving an impression of the tool efficiency and scalability and by considering a case study on bladder tumorigenesis,
focusing on analysis of different phenotypes of the cancer cells.

We aim at both clear and intuitive descriptions, further technical details can be found in the supplementary material.
In addition, scripts and tutorials to reproduce the results are available, integrated in our PyBoolNet software \cite{Klarner2016PyBoolNet}.

\section{Model checking with accepting states}\label{sec:mc}
CTL model checking was originally designed to answer queries by either \emph{yes} or \emph{no} depending on whether the transition system satisfies the query or not.
When model checking is used commercially, i.e., for the verification of software and hardware designs,
it turned out to be useful to analyse the reasons for why a design does not satisfy a specification.
This lead to the development of algorithms that generate counterexamples when queries are false \cite{Clarke2002Counterexamples}.
When model checking is applied to systems biology a similar analysis is desirable.
In particular, the theory which is developed in the subsequent sections requires the capability of computing the exact set of states that satisfies a query, an idea that is more basic then computing counterexamples.
The accepting states of a CTL formula $\varphi$ and STG $(S,\rightarrow)$ are defined by
\[\mathit{Accept}(S,\rightarrow,\varphi) := \{x\in S \mid x\models \varphi\}\]
where $x\models \varphi$ means that $\varphi$ holds in state $x$ of the transition system $(S,\rightarrow)$.
For more details on model checking and CTL, see the supplementary material.
Note that the accepting states are defined without initial states.
By definition, see \cite{Baier2008Principles}, a transition system $(S,\rightarrow)$ satisfies a model checking query with CTL formula $\varphi$
and initial states $I\subseteq S$ iff $I\subseteq\mathit{Accept}(S,\rightarrow,\varphi)$.

For our purposes, the accepting states of one query are often re-used in the CTL formula of a subsequent query.
We therefore need to be able to encode sets of states as Boolean expressions.
Since the state space is already defined in terms of states of Boolean variables,
we may use a disjunctive normal form (DNF) in which each term references a single state of the set.
Clearly, this is not a very efficient encoding and there are usually other, shorter expressions that also represent the set.
We denote by $\mathit{Enc}(X)$ an arbitrary Boolean expression that represents $X\subseteq S$.
Therefore, $\mathit{Enc}(\mathit{Accept}(S,\rightarrow,\varphi))$ is an expression that represents the accepting states of $\varphi$ in $(S,\rightarrow)$.

\section{Basins of attraction}\label{sec:basins}
Basins of attraction describe regions of state space that lead into specific attractors.
Large basins suggest a certain robustness of their attractors against perturbations, see e.g. \cite{davidich2008fissionyeast}, while small basins suggest unlikely attractors.
The study of basins is therefore essential for understanding the likely long-term behaviours of Boolean networks.

The basins of deterministic and non-deterministic STGs are very different
because states in deterministic STGs can always reach exactly one attractor while states in non-deterministic STGs can usually reach several attractors.
It has therefore been suggested, for example in \cite{Tournier2009OperationalInteractions}, to distinguish between weak and strong basins of attraction
where the weak basin contains all states that can reach a given attractor while the strong basin contains all states that can only reach a given attractor.
For an attractor $A\subseteq S$ and an arbitrary state $a\in A$, the weak basin of an attraction is defined by
\[\mathit{WeakBasin}(S,\rightarrow,a) := \{x\in S \mid \exists\pi\in\mathit{Paths}(x):\,a\in\pi\},\]
where $\mathit{Paths}(x)$ is the set of all paths starting from $x$ in the $(S,\rightarrow)$ and $a\in\pi$ means that $a$ appears somewhere along $\pi$.
Note that, as attractors are strongly connected, it is irrelevant which $a\in A$ is chosen.
Biologically, these basins represent the weakest form of commitment to an attractor and the system might never actually reach the given attractor,
either because there are other attractors that are more likely or because the system remains in an intermediate equilibria.

To define the strong basin one may be tempted to require that all paths lead to the attractor.
But this definition is stronger than the notion that a state can only reach $A$.
It excludes states that are on cycles, even if they can, in terms of attractors, only reach $A$.
Nevertheless this type of basin is interesting and since it is cycle-free, we call it $\mathit{CycFreeBasin}(S,\rightarrow,A)$ and define it by
\[\{x\in S \mid \forall\pi\in\mathit{Paths}(x):\,\exists a\in A:\,a\in\pi\}.\]
Note that this definition is in terms of $A$ instead of an arbitrary $a\in A$.
The reason is that in specific cases it makes a difference which $a\in A$ is chosen, e.g., if $A$ is a complex attractor that contains a cycle that does not visit all states of $A$.
These basins form the strictest notion of attraction, the inevitable fall into the given attractor without the possibility of intermediate oscillations or the potential of reaching another attractor.

\begin{figure}[htbp]
\centering
\subfloat[Weak basin]{\includegraphics[width=3.5cm]{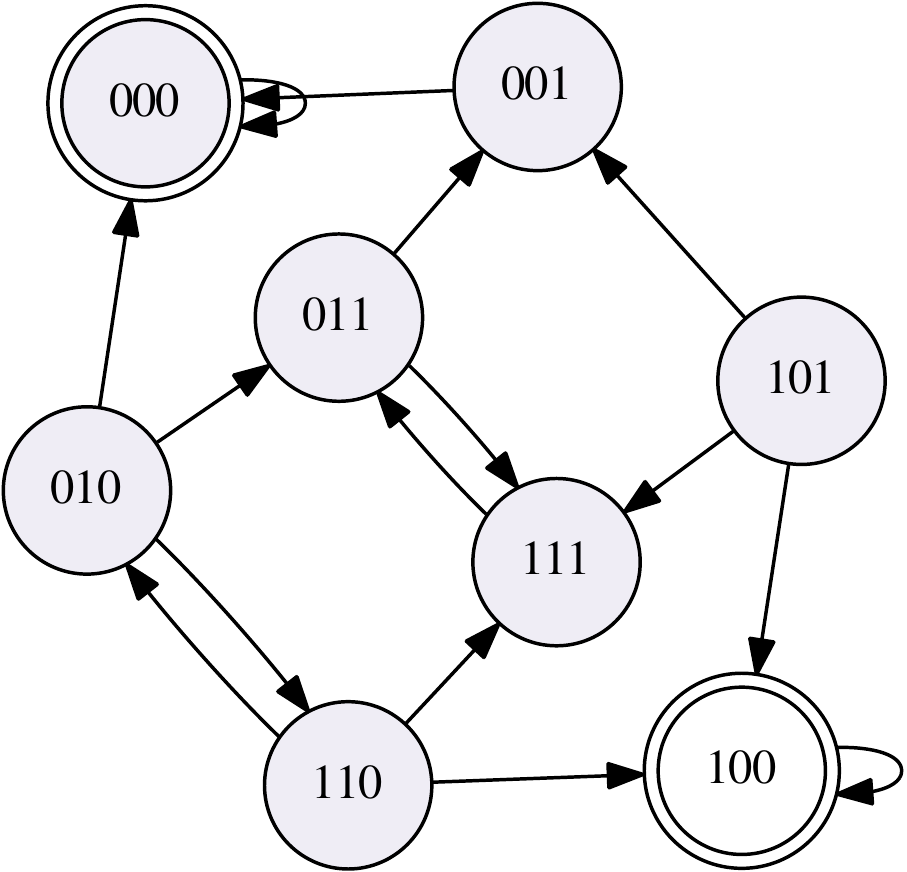}}\hspace{.5cm}
\subfloat[Strong basin]{\includegraphics[width=3.5cm]{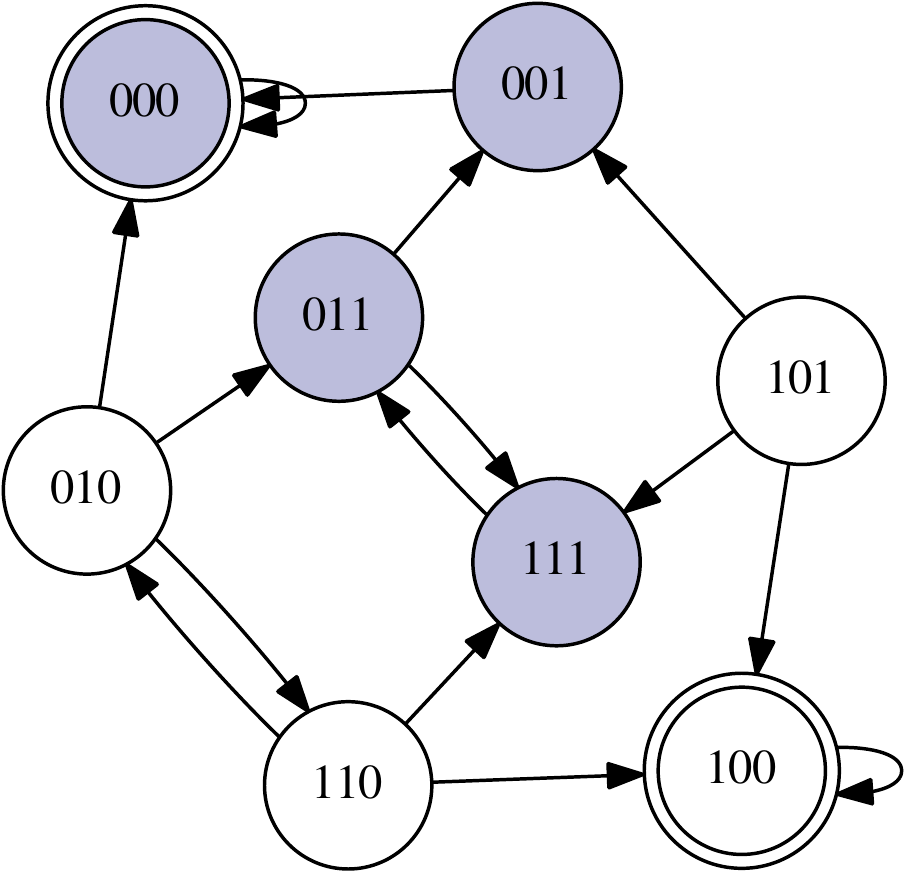}}\hspace{.5cm}
\subfloat[Cycle-free basin]{\includegraphics[width=3.5cm]{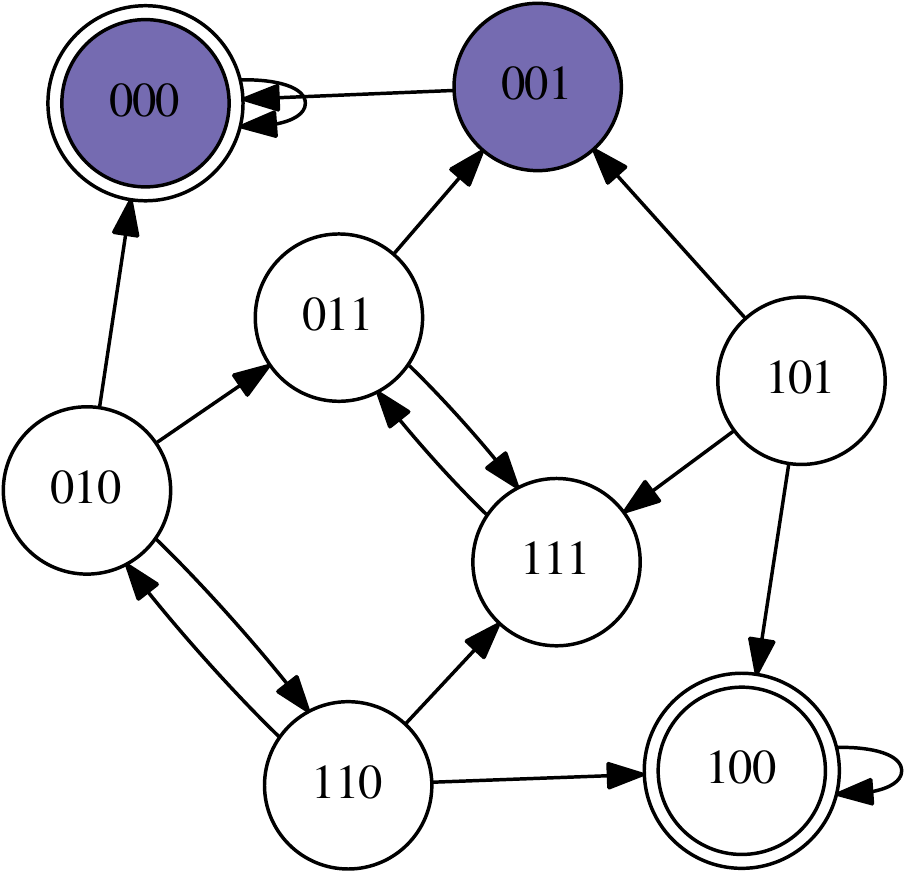}}
\caption{\label{fig:01}
The asynchronous STG of a Boolean network with two steady states highlighted by double circles.
Visualised are the weak, strong and cycle-free basins of $000$. The colours correspond to the ones used in the bar plot of Fig.~\ref{fig:02}~(a).}
\end{figure}

The strong basin in the usual sense of \emph{all states that can only reach $A$} is defined by called $\mathit{StrongBasin}(S,\rightarrow,a)$ and defined by
\[\{x\in S \mid \forall \pi\in\mathit{Paths}(x):\,\forall y\in\pi:\, \exists\pi'\in\mathit{Paths}(y):\, a\in\pi'\}.\]
As the system transitions from a weak to a strong basin it looses the potential of reaching certain other attractors and hence these transitions hint at hidden decision processes within the changing activity patterns.
These three notions are illustrated for a toy transition system in Fig.\ref{fig:01}.
The Boolean equations for the toy transition system were created using iBoolNet \cite{Klarner2018iboolnet}.

To apply symbolic model checking to the computation of weak and strong basins we need CTL queries that capture these definitions.
We only need three operators, namely  the \emph{exists finally} operator $\mathbf{EF}$, the \emph{always finally} operator $\mathbf{AF}$ and the \emph{all globally} operator $\mathbf{AG}$.
These operators can be linked to the path-based definitions of the basins, since a state $x$ satisfies $\mathbf{EF}(M)$ for a state set $M$ if there exists a path starting in $x$ leading to some state in $M$,
similarly $\mathbf{AF}(M)$ encodes a condition for all paths starting in $x$ and the combination $\mathbf{AG}(\mathbf{EF}(M))$ is satisfied if from every state reachable from $x$ there exists a path to a state in $M$.
This last condition would in particular not be true for a state $x$ that can reach an attractor $M$ and in addition a second attractor $M'$, since then the states of $M'$ are reachable from $x$, but $M$ is not reachable from $M'$.
The formal definitions can be found in the supplementary material. These considerations lead us to the following CTL queries
\begin{alignat*}{2}
 \mathit{Accept}(S,\rightarrow,\alpha), \quad\textrm{where}\quad \alpha:=\mathbf{EF}(\mathit{Enc}(a))\\
 \mathit{Accept}(S,\rightarrow,\beta), \quad\textrm{where}\quad \beta:=\mathbf{AF}(\mathit{Enc}(A))\\
 \mathit{Accept}(S,\rightarrow,\gamma), \quad\textrm{where}\quad \gamma:=\mathbf{AG}(\mathbf{EF}(\mathit{Enc}(a)))
\end{alignat*}
for $\mathit{WeakBasin}(a)$, $\mathit{CycFreeBasin}(A)$ and $\mathit{StrongBasin}(a)$, respectively.
Note that the three types of basin satisfy the inclusions
\begin{equation}\label{eq:1}
 A \subseteq \mathit{CycFreeBasin}(A) \subseteq \mathit{StrongBasin}(a) \subseteq \mathit{WeakBasin}(a)
\end{equation}

\begin{figure}[htbp]
\centering
\subfloat[Stacked bar plot]{\includegraphics[width=6cm]{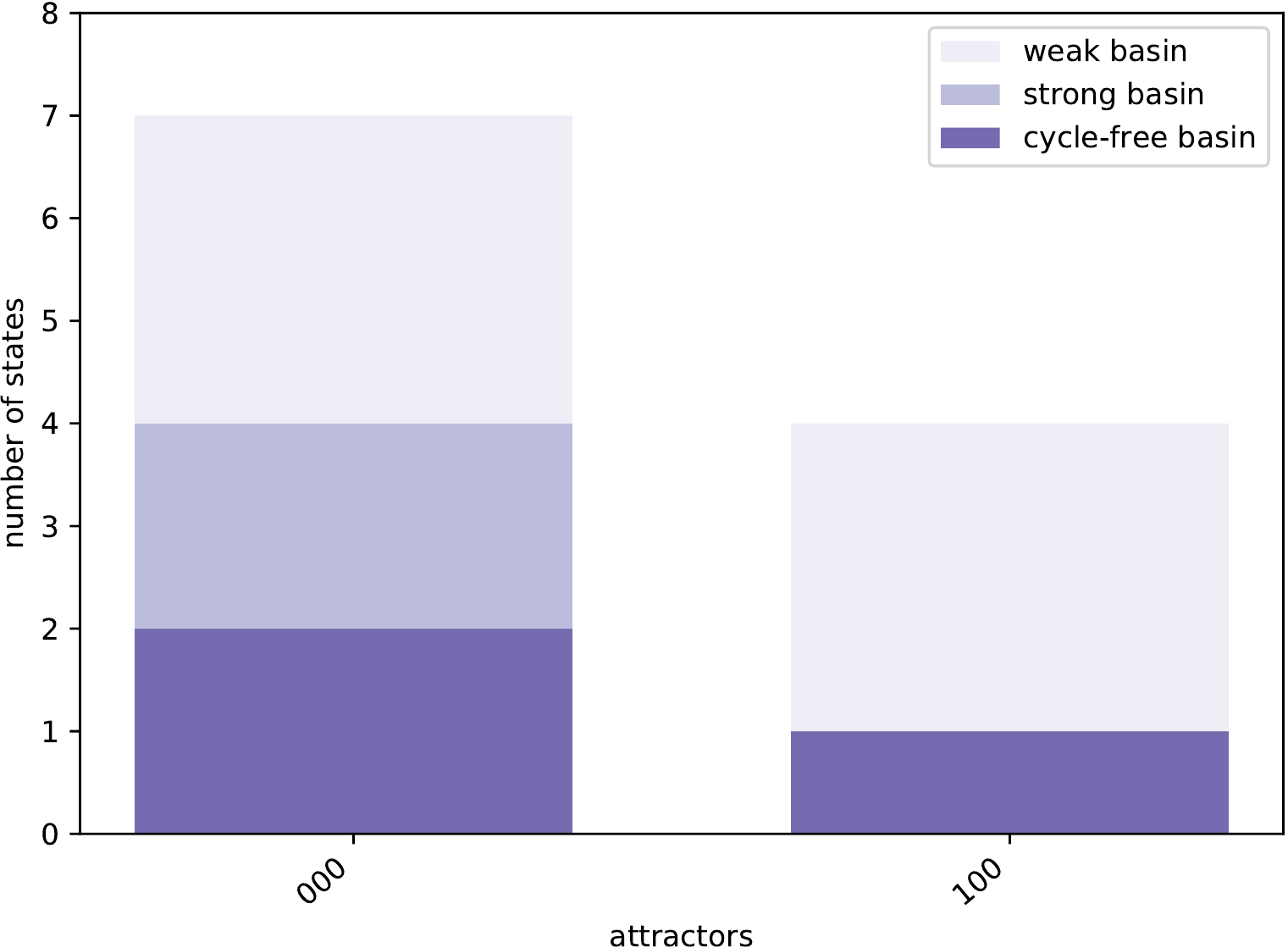}}\hspace{1cm}
\subfloat[Pie chart]{\includegraphics[width=4.5cm]{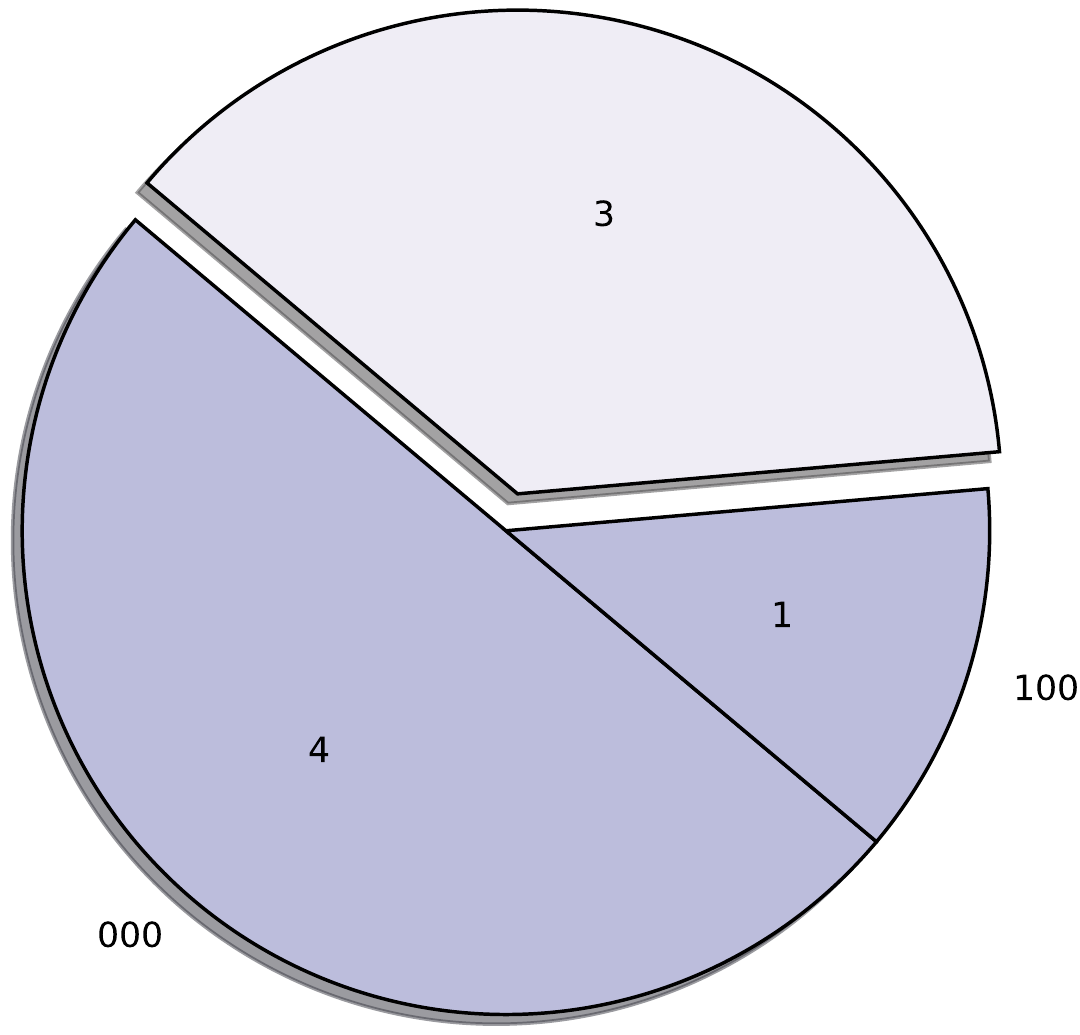}}
\caption{\label{fig:02}
(a) The relative sizes of the basins of attraction as a stacked bar plot that visualised the inclusion relation of Eq.~\ref{eq:1}.
(b) The strong basins of attraction visualised as a pie chart. The light slice represents states that do not belong to a strong basin.
For networks with many states we use percentages rather than absolute number of states.
}
\end{figure}

To visualise the relative sizes of the basins we suggest to use pie charts and stacked bar plots.
The bar plots highlight the inclusion relationship while the pie charts give an impression of the relative sizes of the basins.
These three notions are illustrated for the toy transition system in Fig.\ref{fig:02}.

The notion of basin can be generalised from single to multiple attractors simply by adjusting the state formulas accordingly.
Suppose we are interested in the basins of $k$ attractors $A_1,\dots,A_k$ with representative states $a_i\in A_i$.
With $X:=\{a_1,\dots,a_k\}$ and $Y:=A_1\cup \dots \cup A_k$ the CTL queries for the generalised basins are obtained by using $X$ and $Y$ instead of $a$ and $A$ in the definitions above:
$\mathit{WeakBasin}(X)$, $\mathit{CycFreeBasin}(Y)$ and $\mathit{StrongBasin}(X)$.
Note that while the distributive law $\mathbf{EF}(\phi_1\vee\phi_2)=\mathbf{EF}(\phi_1)\vee\mathbf{EF}(\phi_2)$ for existential CTL operators implies
that the weak basin of several attractors is the union of the individual weak basins:
\[\mathit{WeakBasin}(\{a_i,a_j\}) = \mathit{WeakBasin}(a_i) \cup \mathit{WeakBasin}(a_j)\]
the same does not hold for the strong and cycle-free basins.
Instead, the following inclusions hold:
\begin{alignat*}{2}
 \mathit{CycFreeBasin}(A_i \cup A_j) &\supseteq \mathit{CycFreeBasin}(A_i)\\
 \mathit{StrongBasin}(\{a_i,a_j\})   &\supseteq \mathit{StrongBasin}(\{a_i\})
\end{alignat*}
The reason for the inclusion instead of equality is that a strong basin of two attractors may contain a state that can reach both attractors.
Such a state can not belong to the individual strong basins.
The same example proves that the inclusion may be strict for cycle-free basins.

Note that all definitions above work even if not all attractors are known.
In practice we may therefore compute, for example, the strong basin of a steady state even if we do not know whether there are additional attractors.

\section{Commitment sets}\label{sec:commitment}
A different approach to understanding attractors and regions in state space that lead to them is via commitment sets.
Commitment sets partition the state space into classes that can reach the exact same attractors.
For ease of referencing specific commitment sets, let $(A_1,\dots,A_k)$ be a sequence of $k$ attractors with representatives $(a_1,\dots,a_k)$
and let $I\subseteq\{1,\dots,k\}$ be a non-empty subset of indices.
The states that are committed to the attractors $\{A_i \mid i\in I\}$ are defined by
\[\mathit{ComSet}(I) := \{x\in S \mid (\exists\pi\in\mathit{Paths}(x):\, a_i\in\pi) \Leftrightarrow (i\in I)\}.\]
Commitment sets are the finest resolution of state space with respect to the potential of reaching certain attractors.
An analysis of the similarities of states that belong to the same set is biologically as relevant as understanding the differences between states of different commitment sets.
In particular, it remains an open question to understand the decision processes that are realised by transitions between commitment sets, see also \cite{Berenguier2013HTG}.

To encode this in CTL we observe that a state is committed to the attractors indexed by $I$ if and only if there is a path from $x$ to $a_i$ for every $i\in I$
and $x$ is in the strong basin of $\{A_i \mid i\in I\}$, i.e., can not reach other attractors.
A possible CTL query for $\mathit{ComSet}(I)$ is therefore $\mathit{Accept}(S,\rightarrow,\delta)$, where

\begin{equation*}
\delta:=\bigwedge_{i\in I}\mathit{WeakBasin}(a_i) \wedge \mathit{StrongBasin}(\{a_j \mid j\in I\}). 
\end{equation*}

Note that, as before, commitment sets can be computed without knowing all attractors of the network.
Furthermore, the same queries can be used to compute the basins and commitment sets of subspaces and in particular minimal trap spaces, see \cite{Klarner2015trapspaces},
which are often in a one-to-one correspondence with the attractors of a network, see \cite{Klarner2015Approximations}.

\begin{figure}[htbp]
\centering
\subfloat[Colour-coded commitment]{\includegraphics[width=3.5cm]{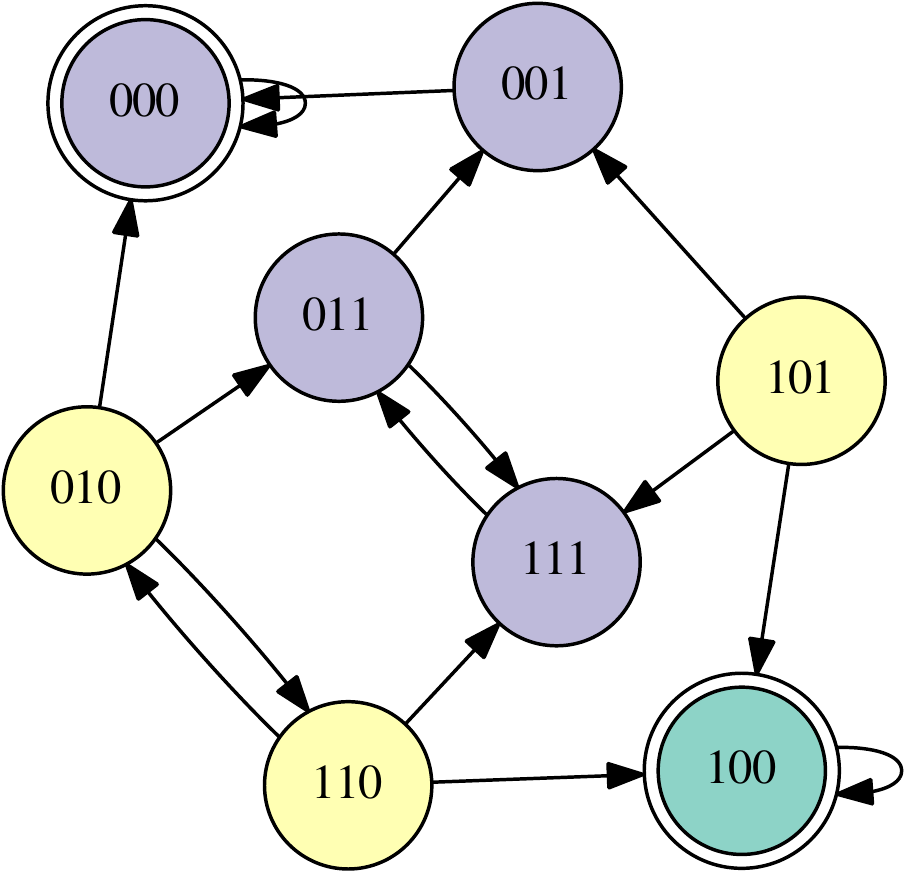}}\hspace{.2cm}
\subfloat[Commitment diagram]{\includegraphics[width=4cm]{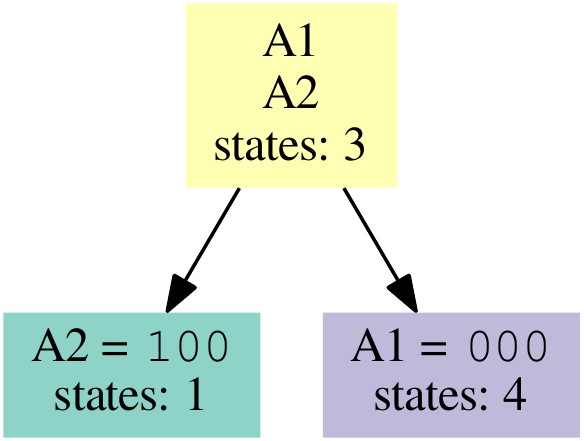}}\hspace{.2cm}
\subfloat[Commitment pie chart]{\includegraphics[width=4cm]{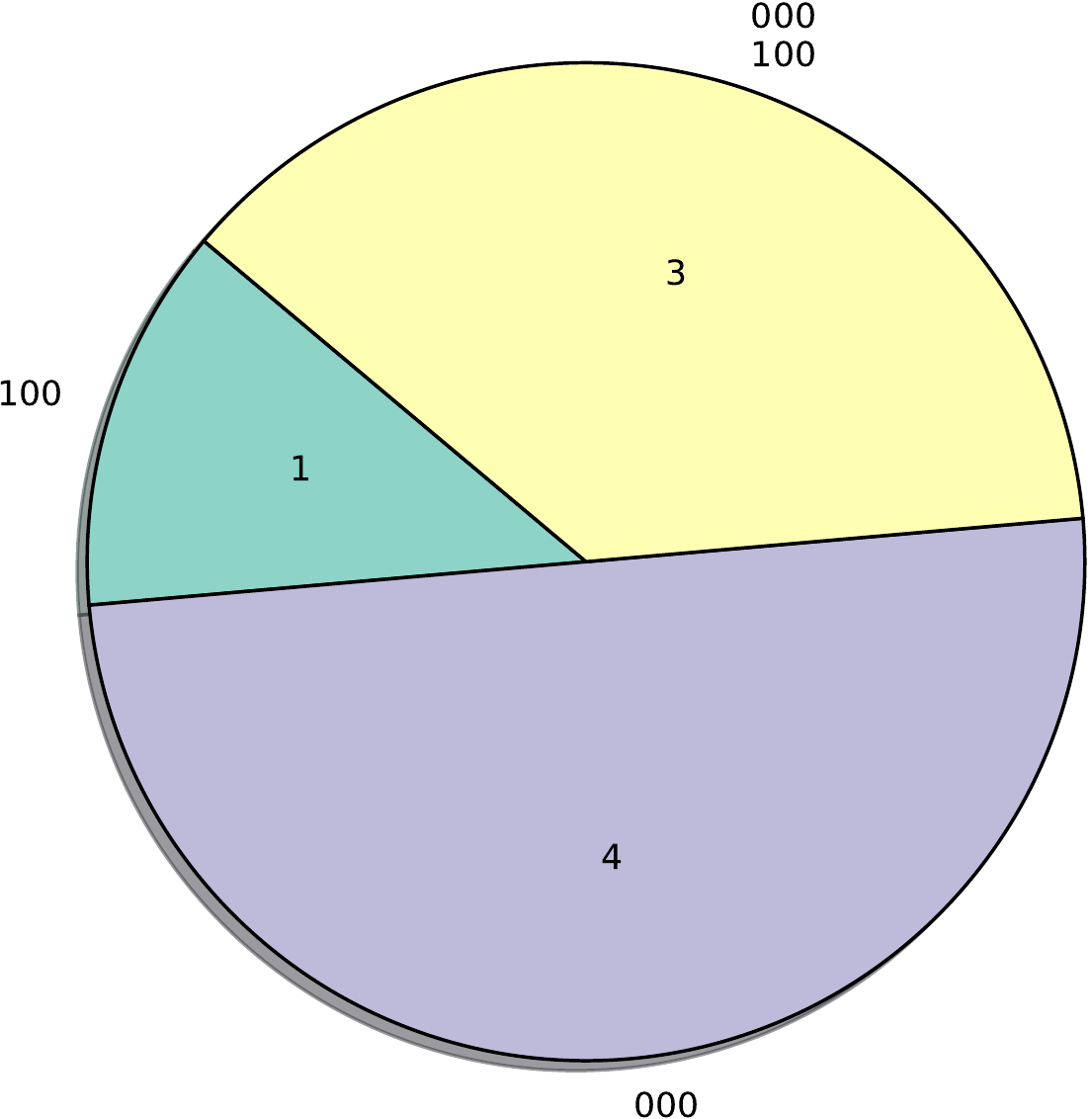}}
\caption{\label{fig:03}
(a) For small STGs we may visualise the commitment sets by using different node colours.
(b) The quotient graph induced by the commitment sets.
(c) The relative sizes of the commitment sets visualised by a pie chart.
}
\end{figure}

For small networks with up to 15 nodes or so we may visualise commitment sets by different node colours.
An example is given in Fig.\ref{fig:03}.
For larger networks we suggest to visualise the commitment sets and transitions between them by a "graph of commitment sets" in which the nodes are commitment sets and there are edges between the nodes if there are transitions between the corresponding commitment sets in the underlying STG.
More precisely, we construct the quotient graph of the STG, induced by the commitment sets.
Recall from basic graph theory that the quotient graph of a digraph $G=(N,\rightarrow)$, induced by a partition $N_1\cup\dots\cup N_m=N$ of its nodes $N$,
is the digraph $Q$ whose nodes are the blocks $N_i$ and there is an edge between
two blocks $N_i\rightarrow N_j$ if there are nodes $n_1\in N_i$ and $n_2\in N_j$ with and edge $n_1\rightarrow n_2$ in $G$.

In the context of commitment sets we call the quotient graph of an STG induced by the commitment sets its commitment diagram.
An example is given in Fig.\ref{fig:03}.

As for the nodes, we compute the edges of the commitment diagram by answering CTL model checking queries, using the "exists next" operator ($\mathbf{EX}$) which is satisfied for a state $x$ if the given formula holds in a state that is reachable from $x$ by a single transition.
Given two commitment sets $\mathit{ComSet}(I)$ and $\mathit{ComSet}(J)$ with $I,J\subset\{1,\dots,k\}$, there is an edge $\mathit{ComSet}(I)\rightarrow \mathit{ComSet}(J)$ if
there are states $x_1\in\mathit{ComSet}(I)$ and $x_2\in\mathit{ComSet}(J)$ such that $x_1\rightarrow x_2$ in the STG.
Using model checking we first compute

\begin{equation*}
\mathit{Accept}(S,\rightarrow,\zeta), \quad\textrm{where}\quad \zeta:=\mathbf{EX}(\mathit{Enc}(\mathit{ComSet}(J))). 
\end{equation*}

and note that the edge $\mathit{ComSet}(I)\rightarrow \mathit{ComSet}(J)$ exists if
\[\mathit{ComSet}(I) \cap \mathit{Accept}(S,\rightarrow,\zeta) \neq\emptyset.\]
Note that an edge can only exist if $I\supsetneq J$ because the set of reachable attractors is decreasing along transitions in an STG.

\section{Phenotypes}\label{sec:phenotypes}
In this section we extend the notion of commitment sets to make it more suited to application.
In experimental studies rarely all components of a system can be measured at the same time, if even at all.
So when using attractors for model validation and prediction an obvious drawback is that the experimental data yields only partial descriptions that might be valid for several attractors.
To be more faithful to this scenario, attractors can be grouped together according to the values of so-called marker components.
The resulting object is called a phenotype.
More specifically, a phenotype is a pattern that specifies the long-term behaviours of a set of given markers with respect to three possible outcomes: \emph{inhibited}=0, \emph{expressed}=1 and \emph{oscillating}=$*$.
A phenotype analysis can help to understand fundamental system properties by simplifying the attractor landscape of the system.
Its main advantage in application, however, is that the measurable components naturally form a set of markers that link the predictions of the Boolean network with observed expression data \cite{WAJaoude2015ThPlasticity,remy2015tumorigenesis}.

We formalise phenotypes as three-valued functions $p:U\rightarrow\{0,1,*\}$ from the set of markers $U\subseteq V$ to the set of possible outcomes $\{0,1,*\}$.
The phenotype $p_A$ of an attractor $A\subseteq S$ is defined by
\[p_A(u)=\begin{cases}0:& \forall a\in A:\, a(u)=0 \\ 1:& \forall a\in A:\, a(u)=1 \\ *:& \textrm{else}\end{cases}\]
Given a set of markers $U$ we say that a phenotype $p$ exists if there is an attractor $A\subseteq S$ such that $p_A=p$.
Note that there may be several attractors with identical phenotypes, $p_A=p_{A'}$, but every attractor has a unique phenotype (for fixed markers).

Analogous to the commitment diagram we define the phenotype diagram by partitioning the state space into subsets that are committed to identical phenotypes.
Let $p_1,\dots,p_k$ be the phenotypes of a given STG and markers.
For each phenotype $p_i$, let $A_{i_1},\dots,A_{i_l}$ be all corresponding attractors with representative states $a_{i_j}\in A_{i_j}$ and denote by $P_i:=\{a_{i_1},\dots,a_{i_l}\}$ the union of the representative states.
We call the commitment set for the phenotypes specified by $I\subseteq\{1,\dots,k\}$ the phenotype set.
It is called $\mathit{PhenoSet}(I)$ and defined by
\[\{x\in S \mid (\exists\pi\in\mathit{Paths}(x):\, \exists a\in P_i:\, a\in\pi) \Leftrightarrow (i\in I)\}.\]
It can be computed using model checking with accepting states $\mathit{Accept}(S,\rightarrow,\eta)$ where

\begin{equation*}
 \eta:=\bigwedge_{i\in I}\mathit{WeakBasin}(P_i) \wedge \mathit{StrongBasin}(\bigcup_{i\in I} P_i). 
\end{equation*}
Note the similarity with the query for commitment sets.
The only difference is that we use sets of representatives $P_i$ instead of individual representatives $a_i$.
As before, the existence of edges $\mathit{PhenoSet}(I)\rightarrow \mathit{PhenoSet}(J)$ between phenotype sets is determined via the quotient graph relation.

\section{NuSMV with Accepting States}\label{sec:nusmva}
To compute the basin, commitment and phenotype diagrams, we developed an extension of \emph{NuSMV} that is capable of returning the accepting states of a transition system and CTL formula.
It is based on \emph{NuSMV} $2.6$, called \emph{NuSMV-a} ("a" for "accepting states"), available from \cite{Klarner2016NuSMVa} and fully integrated in \emph{PyBoolNet} starting from version $2.0$ \cite{Klarner2016PyBoolNet}.
In order to decide whether a CTL formula $\varphi$ is true or false,
\emph{NuSMV} already computes the accepting states in the form of a BDD that represents $\mathit{Accept}(S,\rightarrow,\varphi)$.
We only had to make this information accessible and implement a suitable format in which to return the accepting states.
We decided to use an output format called \emph{factored form} that is implemented in the function \texttt{Cudd\_DumpFactoredForm} (\emph{CUDD} $3.1.4$).
Although the resulting Boolean expressions become quickly unreadable for humans they are in \emph{NuSMV} syntax and can therefore be used as sub-formulas of subsequent CTL queries.
To obtain minimised expressions that may be easier for humans to read we extended PyBoolNet with a module called \texttt{BooleanLogic} and
specifically the function \texttt{minimize\_espresso(..)} that uses the program \emph{espresso} for minimising Boolean expressions \cite{brayton1984espresso}.

As of version $2.0$, \emph{PyBoolNet} comes with \emph{NuSMV-a} and the function
\[\texttt{check\_primes\_with\_acceptingstates(..)}\] of the package \texttt{ModelChecking} returns
a \emph{Python} dictionary that contains all formulas and cardinalities as strings and integers, respectively.
The PyBoolNet modules for computing the various basins, the commitment sets and the phenotypes all require this new model checking function.
For tutorials and examples consult the PyBoolNet manual and project homepage.

\section{Efficiency and Scalability}\label{sec:antelope}
This section presents a table of running times for computing the accepting states of four Boolean networks of increasing size:
\texttt{tournier\_apoptosis} \cite{Tournier2009OperationalInteractions} with $n=12$,
\texttt{dahlhaus\_neuroplastoma} \cite{dahlhaus2016neuroplastoma} with $n=24$,
\texttt{grieco\_mapk} \cite{Grieco2013MAPK} with $n=53$ and
\texttt{jaoude\_thdiff} \cite{WAJaoude2015ThPlasticity} with $n=103$.
The networks are available in the PyBoolNet model repository.
We generated three CTL queries of increasing complexity for each network: a weak basin query, a strong basin query and a commitment set query.
For state formulas we compute two attractor states $a_i\in A_i$ for $i=1,2$ for each network using minimal trap spaces and used a conjunctive encoding.
We compared the running times of NuSMV-a with \emph{Antelope} (Analysis of Networks through Temporal-Logic Specifications) the only other software we found for computing accepting states of CTL queries.
\emph{Antelope} is a model checking software specifically designed for biological interaction networks, see \cite{Arellano2011Antelope}.
Its model checking algorithms are based on "hybrid CTL", an extension of standard CTL that allows, for example, the specification of cycles of a given length.
We used a 64-bit Ubuntu 16.04 LTS desktop pc with 7.7 GiB memory and 4 CPUs with 2.7 GHz and the NuSMV-a command \texttt{-dcx -dynamic -df -coi -a print}.
The results in Tab.~\ref{tab:01} show that NuSMV can handle networks with up to around 50 variables in reasonable times in the order of minutes.
Antelope can answer queries in reasonable time for about 20 variables.

\begin{table*}[htbp]
\centering
\begin{tabular}{clcc}
\hline
\multicolumn{1}{l}{}                                                                                                    &                          & \multicolumn{1}{l}{NuSMV-a} & \multicolumn{1}{l}{Antelope} \\ \hline
\rowcolor[HTML]{EFEFEF} 
\cellcolor[HTML]{EFEFEF}                                                                                                & \texttt{weak\_basin}     & 1s                          & 1s                           \\
\rowcolor[HTML]{EFEFEF} 
\cellcolor[HTML]{EFEFEF}                                                                                                & \texttt{strong\_basin}   & 1s                          & 1s                           \\
\rowcolor[HTML]{EFEFEF} 
\multirow{-3}{*}{\cellcolor[HTML]{EFEFEF}\begin{tabular}[c]{@{}c@{}}\texttt{tournier\_apoptosis}\\ $n=12$\end{tabular}} & \texttt{commitment\_set} & 1s                          & 1s                           \\
                                                                                                                        & \texttt{weak\_basin}     & 1s                          & 19m6s                        \\
                                                                                                                        & \texttt{strong\_basin}   & 1s                          & 18m3s                        \\
\multirow{-3}{*}{\begin{tabular}[c]{@{}c@{}}\texttt{dahlhaus\_neuroplastoma}\\ $n=24$\end{tabular}}                     & \texttt{commitment\_set} & 1s                          & 20m31s                       \\
\rowcolor[HTML]{EFEFEF} 
\cellcolor[HTML]{EFEFEF}                                                                                                & \texttt{weak\_basin}     & 1m45                        & $>$12h                       \\
\rowcolor[HTML]{EFEFEF} 
\cellcolor[HTML]{EFEFEF}                                                                                                & \texttt{strong\_basin}   & 2m15                        & $>$12h                       \\
\rowcolor[HTML]{EFEFEF} 
\multirow{-3}{*}{\cellcolor[HTML]{EFEFEF}\begin{tabular}[c]{@{}c@{}}\texttt{grieco\_mapk}\\ $n=53$\end{tabular}}        & \texttt{commitment\_set} & 3m51                        & $>$12h                       \\
                                                                                                                        & \texttt{weak\_basin}     & $>$12h                      & $>$12h                       \\
                                                                                                                        & \texttt{strong\_basin}   & $>$12h                      & $>$12h                       \\
\multirow{-3}{*}{\begin{tabular}[c]{@{}c@{}}\texttt{jaoude\_thdiff}\\ $n=103$\end{tabular}}                             & \texttt{commitment\_set} & $>$12h                      & $>$12h                      
\end{tabular}
\caption{
The running times for answering CTL queries with accepting states for four Boolean networks of increasing size.
See main text for details.
}\label{tab:01}
\end{table*}

\section{Case Study: Bladder Tumorigenesis}\label{sec:tumorigenesis}
In this case study we investigate the bladder tumorigenesis model of \cite{remy2015tumorigenesis} and compare the results the trajectory-based studies that are also available in \cite{remy2015tumorigenesis}.
The tumorigenesis model has four inputs which represent the cellular context of the network in terms of
DNA damage, the growth factors EGFR and FGFR3 and growth inhibitors. Three output components, Growth
Arrest, Proliferation and Apoptosis determine the response of the cell to a given context.
The model was constructed using the software GINsim \cite{Naldi2009Ginsim}.
There are five multi-valued components: ATM, Apoptosis, CHEK1\_2, E2F1 and E2F3 each with the activity levels 0=low, 1=medium and 2=high.
Since we use NuSMV-a as a plugin of PyBoolNet which is designed for Boolean networks we worked with the equivalent Booleanised version in which multi-valued components x are replaced by two Boolean components x\_medium and x\_high.
The update functions are obtained by the so-called van Ham mapping which has been shown to result in an equivalent transition system \cite{Didier2011Mapping}.
PyBoolNet detects variables that represent the same multi-valued variable and removes inadmissable combinations of values from the initial states of model checking queries.
Meaningless states, e.g., x\_medium=0 and x\_high=1, that are artefacts of the Booleanisation do therefore not contribute to the basins, commitment sets and phenotype sets.
The Booleanisation was done using GINsim and the Boolean model is available from the PyBoolNet model repository under the name remy\_rumorigenesis.

The analysis of the wild-type multi-valued model consists of a list of the observed attractors and their phenotypes, see Table~S2 of the supplementary materials of \cite{remy2015tumorigenesis}, here referred to as Table~S2,
and pie charts that visualise the reachability of each phenotype under different conditions, e.g. the wildtype but also various knock-out mutants.
The reachabilities of phenotypes are defined as the percentage of simulations (random walks on the asynchronous state transition graph) that exhibit a certain phenotype.
The authors ran 10,000 simulations for each of the 16 input configurations where the remaining variables are randomly initialised.

We first confirmed that the attractors, 20 steady states and 5 cyclic attractors, of Table~S2 are identical to the ones we find using the minimal trap space method \cite{Klarner2015trapspaces, Klarner2015Approximations} for our Boolean model.
The 25 attractors exhibit four phenotypes for the markers Proliferation, Apoptosis\_medium, Apoptosis\_high and Growth\_arrest.
The phenotypes are listed in Tab.~\ref{tab:02}.

\begin{figure}[htbp]
  \centering
  \includegraphics[width=\linewidth]{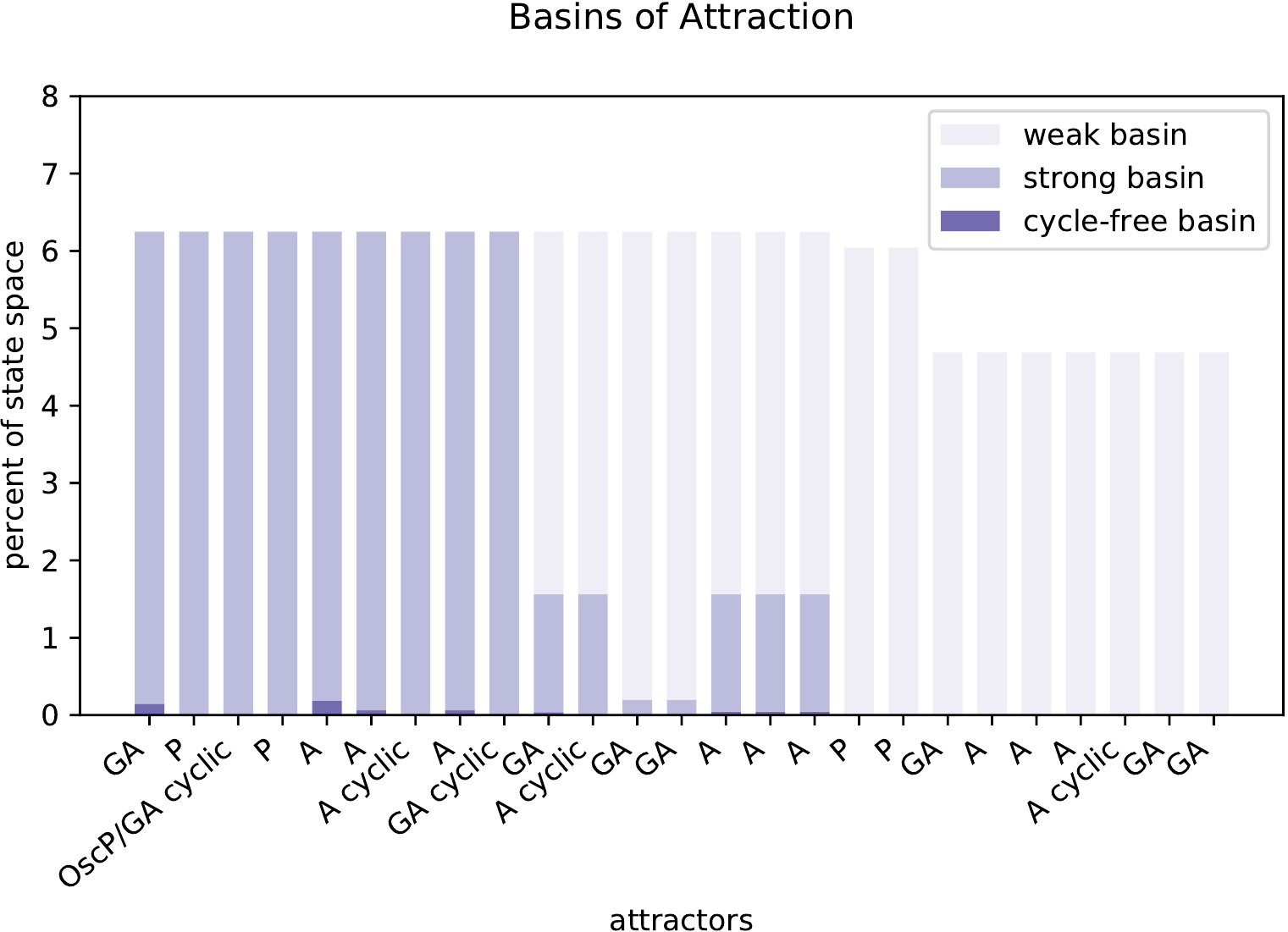}
  \caption{\label{fig:06}
  The size of the weak, strong and cycle-free basins of the 25 attractors.
  Each bar is labelled by the phenotype pattern of the corresponding attractor. Note that basins whose size is too small compared to the state space size cannot be detected in the plot.
  }
\end{figure}

The basins of attraction are visualised as a bar plot in Fig.~\ref{fig:06}.
Since the model has four inputs, its state space is divided into 16 disconnected regions each making up $1/16=6.25\%$ of all states.
The maximal height of a bar is therefore at $6.25$.
It can be deduced that the first 9 attractors are the only ones in the state space region defined by their respective input combination because the whole region belongs to their strong basins.
The remaining 16 attractors share their input combination with another attractor and the majority of states belong to more than one basin.

\begin{figure*}[htbp]
  \centering
  \subfloat[]	{\includegraphics[width=\linewidth]{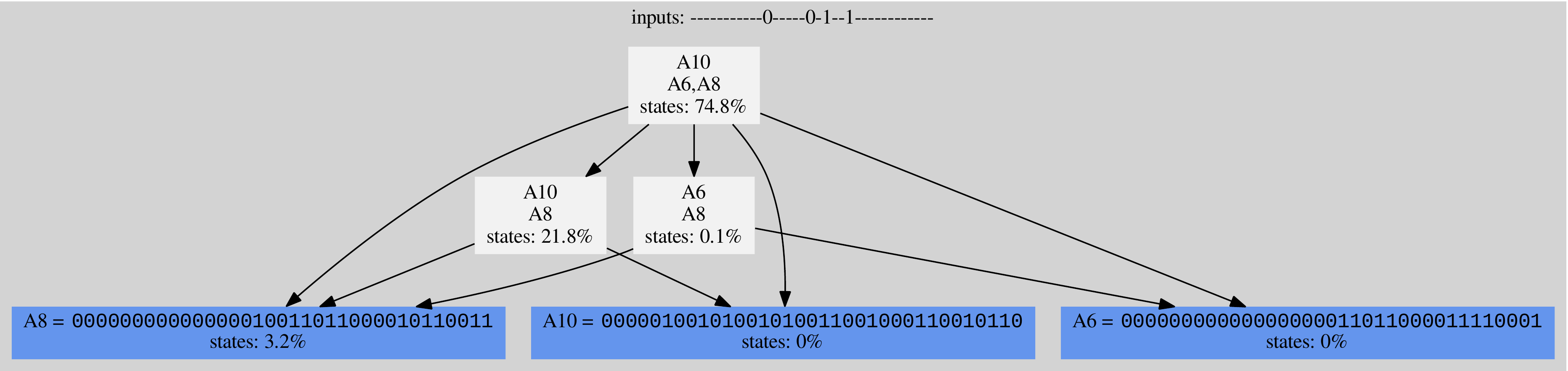}}
  
  \subfloat[]	{\includegraphics[width=0.6666\linewidth]{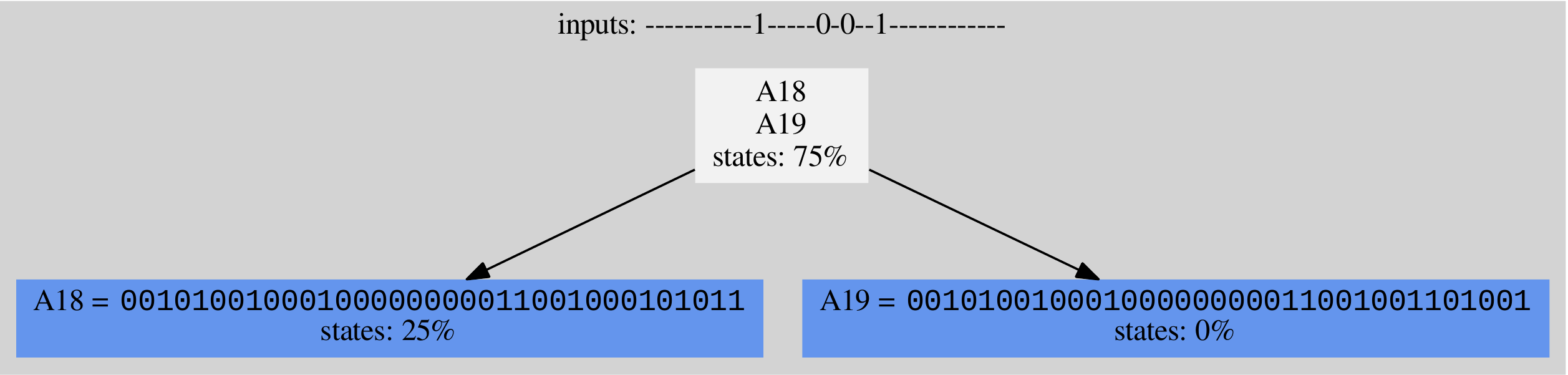}}
  \subfloat[]	{\includegraphics[width=0.3333\linewidth]{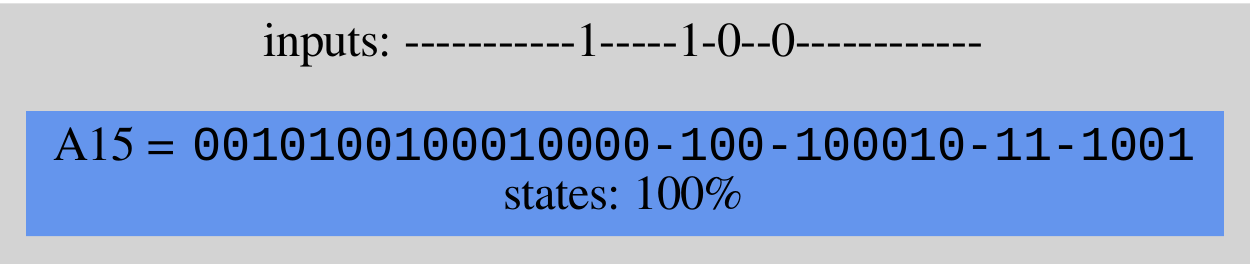}}\\
  
  \caption{\label{fig:07}
  Three of the 16 parts of the commitment diagram. Each of the remaining 13 parts is isomorphic to one of these, including the size of the nodes.
  }
\end{figure*}

The commitment diagram also reflects the division into 16 disconnected regions.
It is divided into 16 disconnected diagrams of three different shapes where the shape depends only on how many attractors there are in the respective region.
An example of each is given in Fig.~\ref{fig:07}.
If an input combination contains more than one attractor, Fig.~\ref{fig:07}~(a) and (b), then the majority of states, about 75\%, can reach all attractors.

\begin{figure*}[htbp]
  \centering
  \subfloat[]	{\includegraphics[width=6.1cm]{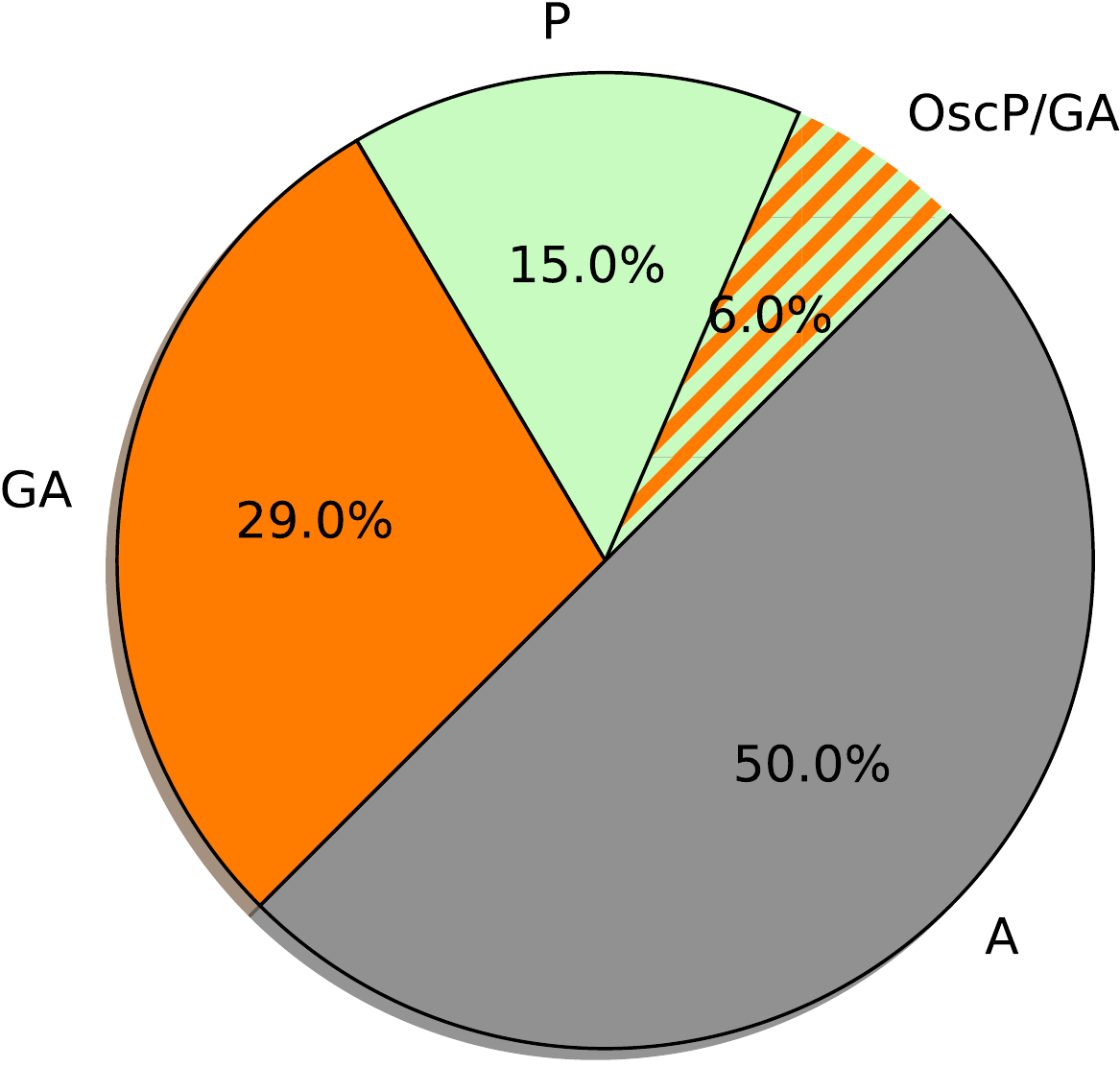}}
  \subfloat[]	{\includegraphics[width=6cm]{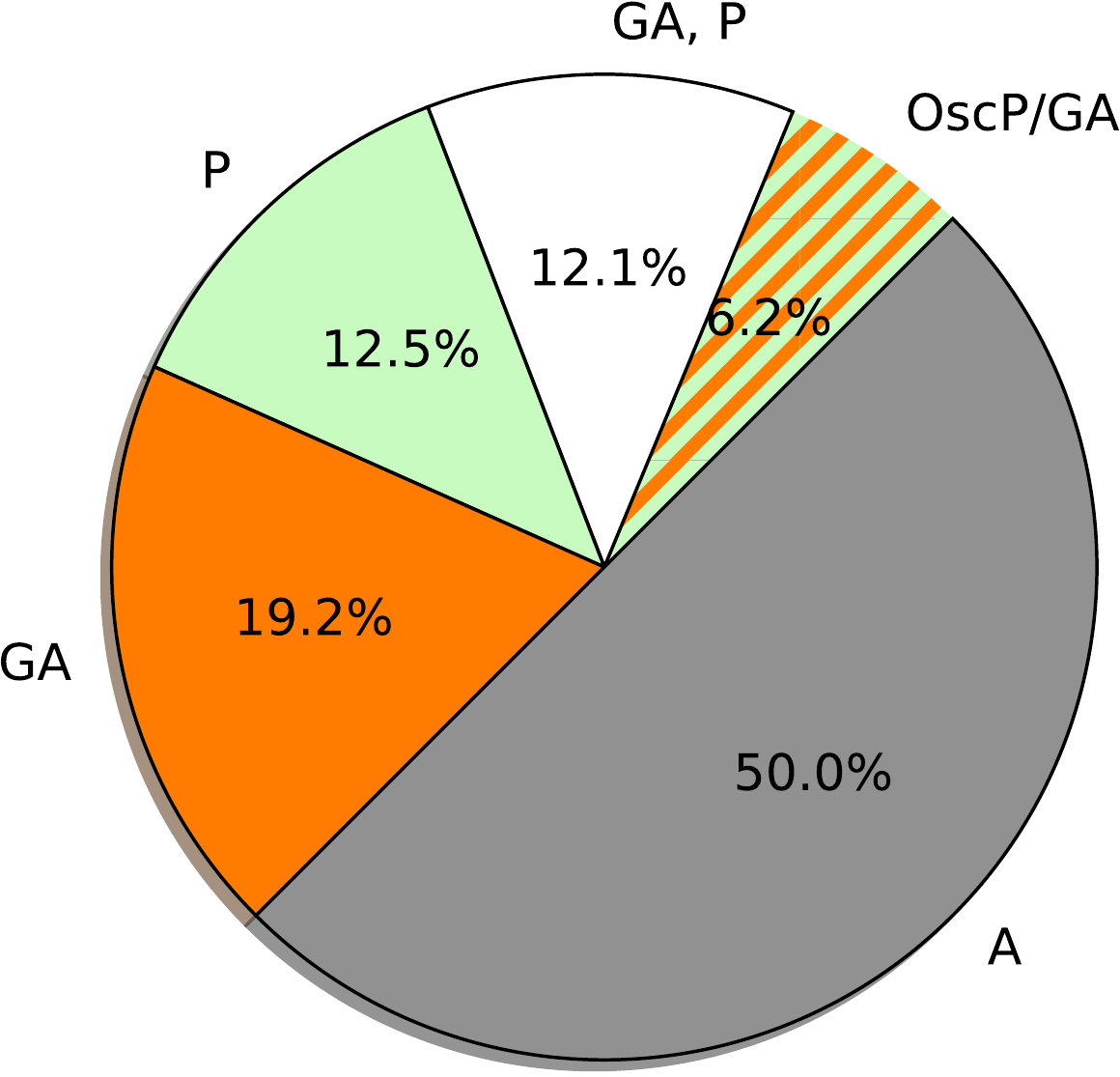}}\\
  
  \caption{\label{fig:05}
  (a) The trajectory-based likelihood of reaching one of the four phenotypes reproduced from \cite{remy2015tumorigenesis}.
  (b) The fractions of state space that are committed to identical sets of phenotypes.
  }
\end{figure*}

We then compared the phenotype reachability analysis of \cite{remy2015tumorigenesis} with our notion of phenotype commitment.
The reachability analysis computes the percentage of simulations that exhibit a given phenotype.
The results are visualised as a pie chart with one slice per phenotype, see Fig.~\ref{fig:05}~(a) (taken from \cite{remy2015tumorigenesis}).
Our state-based analysis, which uses model checking with accepting states, finds that there are five different phenotype commitment sets:
one for each phenotype but also one additional slice that represents states in the basin of several phenotypes, namely   \texttt{GA} and \texttt P, see Fig.~\ref{fig:05}~(b).
To make the comparison between the trajectory-based analysis of \cite{remy2015tumorigenesis} and our state-based analysis easier we also visualise the commitment sets using a pie chart and use the same colours.
Looking at both charts suggests that the majority of the states not yet committed to a single attractor (namely $10/12\%$) lead to \texttt{GA} but neither chart alone holds the full information of the probability of reaching the phenotypes and the sizes of the commitment sets.
The two types of analysis can be seen as complementing each other.

\begin{table}[htbp]
\centering
\begin{tabular}{ccccccc}
\hline
                                                              & \multicolumn{1}{l}{$m_1$} & \multicolumn{1}{l}{$m_2$} & \multicolumn{1}{l}{$m_3$} & \multicolumn{1}{l}{$m_4$} & \multicolumn{1}{l}{steady} & \multicolumn{1}{l}{cyclic} \\ \hline
\cellcolor[HTML]{FF7C00}\verb+GA+                             & \multicolumn{1}{c|}{0}    & 0                         & \multicolumn{1}{c|}{0}    & 1                         & 7                          & 1                          \\
\cellcolor[HTML]{C8FBC0}\verb+P+                              & \multicolumn{1}{c|}{1}    & 0                         & \multicolumn{1}{c|}{0}    & 0                         & 4                          & 0                          \\
\cellcolor[HTML]{919191}\verb+A+                              & \multicolumn{1}{c|}{0}    & 1                         & \multicolumn{1}{c|}{0}    & 1                         & 9                          & 3                          \\
\cellcolor[HTML]{C8FBC0}{\color[HTML]{FF7C00} \verb+OscP/GA+} & \multicolumn{1}{c|}{*}    & 0                         & \multicolumn{1}{c|}{0}    & *                         & 0                          & 1                         
\end{tabular}
\caption{
The four phenotypes of the wildtype model with markers $m_1=$ Proliferation, $m_2=$ Apoptosis\_medium, $m_3=$ Apoptosis\_high and $m_4=$ Growth\_arrest.
The colour coding is consistent with Fig.~\ref{fig:05}.
The last two columns record the number of steady states and cyclic attractors of each phenotype.
Note that a phenotype may contain cyclic attractors even though all markers are steady.
}\label{tab:02}
\end{table}

\section{Discussion}\label{sec:discussion}

This paper puts the notions of a basin of attraction for non-deterministic transition systems on firm mathematical ground.
We propose definitions of varying strictness and focus that allow to partition state space in a way providing insight in the sequences of both reversible and non-reversible decisions ultimately leading the system into an attractor.
For each introduced set and its corresponding reachability properties, we developed means of visualisation, either in the form of a diagram or a pie chart.

On the computational side we define what CTL model checking with accepting states is and offer an extension to NuSMV that is capable of returning the accepting states of a query.
This extension, NuSMV-a, as well as the computation and visualisation tools for all notions introduced here are now part of PyBoolNet which offers a high level interface to model checking Boolean networks with many functions for model and query generation already implemented.
The comparison of NuSMV-a with Antelope suggests that computing accepting states with NuSMV-a is competitive and networks with more than 50 nodes can be handled on a scale of minutes.

Using the case study, we illustrated the various diagrams for a molecular network and also to confirm the correctness of the results in comparison with the study of Remy et al., e.g., concerning the calculated attractors and phenotypes.
We showed that our approach yields a useful complementary view on assessing the likelihood of reaching particular attractors and the corresponding variability inherent in the system.

We believe based on our results that in the context of non-deterministic transition systems a range of notions concerning attractor reachability are of interest.
The cycle-free basins springing from the natural condition that all paths originating in a state should lead to the same attractor show that further distinctions concerning the asymptotic behavior of trajectories are possible and useful.
The question in how far non-terminal cycles in the state transition graph can play a similar role to attractors lead to questions about the stability of such cycles.
Progress concerning this issue in turn could potentially be highly meaningful for application, since it may allow to uncover robustly maintainable system behaviors in the models that are overlooked when only focusing on the attractors.
A possible way to address this question mathematically and computationally may be in combining our state based analysis of commitment sets with trajectory-based analysis as we have outlined in the case study for the reachability analysis of phenotypes.

In a similar step to higher resolution, the commitment sets can be further refined by considering not only which attractors are reachable but also which other commitment sets, i.e.,
two states may both be committed to the same attractors but they traverse different commitment sets on their way to the attractors.
It is certainly of mathematical interest to find the rule allowing to link the properties of the resulting graphs on the different resolution levels, which in turn links to the most exciting open question that results from this research is the question of characterising the control points and decision processes that occur at the border between commitment sets.
Which changes in activity of which variables are responsible for the loss of a previously attainable attractor?
What are minimal interventions that restore the reachability?
Questions along these lines are crucial for understanding the dynamics of qualitative models and in practice for controlling cell differentiation and cell reprogramming.

\section*{Acknowledgment}

We thank the participants of the discussion about basins of attraction at Luminy University in Marseille in 2016:
Elisabeth R\'emy, Claudine Chaouiya, Laurence Calzone, Denis Thieffry, Aur\'elien Naldi, Pedro Monteiro and Laurent Tournier.
The work was partially funded by the German Federal Ministry of Education and Research (BMBF), grant no. 0316195.

\bibliographystyle{unsrt}
\bibliography{klarnerbib}

\end{document}